\documentclass[12pt]{amsart}
\usepackage{amsmath}

\newtheorem{thm}{Theorem}
\newtheorem{lem}[thm]{Lemma}
\newtheorem{prop}[thm]{Proposition}
\newtheorem{defin}[thm]{Definition}

\theoremstyle{definition}

\theoremstyle{remark}

\newcommand{\Z}{{\mathbf Z}}
\newcommand{\N}{{\mathbf N}}

\newcommand{\R}{{\mathbf R}}

\newcommand{\tr}{\operatorname{trace}}
\newcommand{\HH}{\mathbf H} 
\newcommand{\TT}{\mathbf T} 

\newcommand{\norm}[1]{\lVert #1\rVert}
\newcommand{\eps}{{\epsilon}}

\begin{document}

\title
[Stable mixing and quasi-morphisms]
{Stable mixing for cat maps and quasi-morphisms
of the modular group}
\author{Leonid Polterovich and Zeev Rudnick}
\address{ School of Mathematical Sciences,
Tel Aviv University, Tel Aviv 69978, Israel
({\tt polterov@math.tau.ac.il})}
\address{School of Mathematical Sciences,
Tel Aviv University, Tel Aviv 69978, Israel
({\tt rudnick@math.tau.ac.il})}

\date{July 14, 2003}
\thanks{Supported by THE ISRAEL SCIENCE FOUNDATION
founded by the Israel Academy of Sciences and Humanities.
}
\maketitle

\section{Introduction}

It is well-known that the action of a hyperbolic element
(``cat map'')  of the modular group $ SL(2,\Z)$ on the torus
$\TT = \R^2/\Z^2$ has
strong  chaotic dynamical properties such as
mixing and exponential decay of correlations.
In this note we study stability of this behaviour with respect to a
class of perturbations called ``kicked systems''  introduced in
\cite{kicks1}: Let $h\in SL(2,\Z)$
be hyperbolic, and assume that the
system $\{h^n\}, \; n \in \N$ is influenced by a sequence of kicks
$\Phi=\{\phi_1,\phi_2,\dots\} \subset SL(2,\Z)$ which arrive
periodically with
period $t \in \N$. The evolution of the kicked system is given by the
sequence of maps
\begin{equation}\label{kicked system}
f^{(n)}(t) =
\phi_{n}h^t\dots \phi_2 h^t\phi_1 h^t.
\end{equation}
We shall always assume that the sequence $\tr (\phi_i)$
is bounded.

The kicked system defined in \eqref{kicked system}, is a particular
case of the more general notion of a sequential dynamical system. Any
sequence of maps $f^{(n)}$ gives rise to  dynamics on the torus, the
integer $n$ playing the role of discrete time, and the evolution of
an observed quantity, that is a function $F\in L^2(\TT)$ is given by
$F\circ f^{(n)}$.
Given two observables $F_1,F_2\in L^2(\TT)$, the
time correlation function is
\begin{multline}\label{correlation function}
C(F_1,F_2;f^{(n)}) = \\
\int_{\TT} F_1(f^{(n)}(x)) F_2(x) dx - \int_{\TT}F_1(x)dx
\int_{\TT}F_2(x) dx \;.
\end{multline}
The system is {\em mixing} if $C(F_1,F_2;f^{(n)}) \to 0$ as
$n\to \infty$ for all $F_1,F_2\in L^2(\TT)$. For instance,
the sequence  $\{h^n\}, n \in N$ is mixing if and only if
$h$ is hyperbolic.

\begin{defin}
An  element $h\in SL(2,\Z)$ is stably mixing if for every
sequence of kicks $\Phi=\{\phi_i\}$ with bounded traces there exists
$t_0=t_0(\Phi)$ so that for every $t>t_0$, the sequential system
$f^{(n)}(t) =\phi_{n}h^t\dots \phi_2 h^t\phi_1 h^t$ is mixing.
\end{defin}

As it was noticed in
 \cite{kicks1} if $h$ is conjugate to its inverse in $SL(2,\Z)$, it
is {\em not} stably mixing, since if $ghg^{-1}=h^{-1}$ and we take as
the sequence of kicks $g^{-1},g,g^{-1},g,\dots$ then the evolution of
the kicked system $f^{(n)}(t)$ is $2$-periodic: $f^{(n)}(t)$ is either
$g^{-1}h^t$ or $1$ depending if $n$ is odd or even.

In this paper we give a complete characterization for stable mixing
for the 2-torus:
\begin{thm}\label{thm: mixing}
A hyperbolic element $h\in SL(2,\Z)$ is stably mixing if and only if
it is not conjugate to its inverse in $SL(2,\Z)$.
\end{thm}
Moreover we give a quantitative measure of mixing by showing
(Theorem~\ref{thm: decay of correlations} below)
that the correlations decay
exponentially for H\"older observables, as well as establish
mixing for a more general class of perturbations (see \S 4).
In the last section we study ``Lyapunov exponents" of kicked systems.

The proof of Theorem \ref{thm: mixing}
  is a combination of basic harmonic
analysis and geometric group theory.
The essential notion that we
need is that of a {\em quasi-morphism} (see e.g. \cite{quas}): Given
a group $G$,
and a function $r:G\to \R$, we say that $r$ is a {\em quasi-morphism}
if its {\it derivative} $dr:G \times G \to \R$ defined by
$$dr(g_1,g_2) = r(g_1g_2)-r(g_1)-r(g_2)$$
is a {\em bounded} function.
A {\em homogeneous} quasi-morphism also satisfies $r(g^n)=nr(g)$
for all $n \in \Z$.

In our previous work \cite[1.5.D]{kicks1}
(see also \cite{EF},\cite{BF}), we showed that given a
hyperbolic element $h$ of $SL(2,\Z)$, not conjugate to its inverse,
there is a homogeneous quasi-morphism of $SL(2,\Z)$ which
does not vanish on $h$.
Crucial to our purpose is the following refinement of that result:
\begin{prop}\label{prop: quasimorphism}
Let $h\in SL(2,\Z)$ be hyperbolic, not conjugate to its inverse. Then
there exists a homogeneous quasi-morphism $r$ of $SL(2,\Z)$ such that
$r(h)=1$ and which vanishes on all parabolic elements.
\end{prop}
This is a special case of the following ``separation theorem'':
\begin{thm}\label{separation}
Let $G$ be a discrete subgroup of $PSL(2,\R)$ and $h\in G$ a primitive
element of infinite order, not conjugate to its
inverse in $G$. Then given any
finite set of elements $g_1,\dots, g_M\in G$, not conjugate to a power
of $h$, there is a homogeneous
quasi-morphism $r$ so that $r(h)\neq 0$, while $r(g_1)=\dots
=r(g_M)=0$.
\end{thm}
The existence of such a rich set of quasi-morphisms is a remarkable
property of discrete groups. For instance, $\Z^n$ fails to satisfy it
for $n\geq 2$ (all homogeneous quasi-morphisms are homomorphisms in
that case).
Moreover, for $SL(n,\Z)$, $n\geq 3$ there are {\em no} nontrivial
quasi-morphisms!  This is a consequence of the bounded generation
property - every matrix in $SL(n,\Z)$ is a product of a bounded number
of elementary matrices
\footnote
{An elementary matrix has 1's on the diagonal, while
all but one off-diagonal entries vanish.}
if $n\geq 3$ \cite{CK1,CK2} and for $n\geq 3$
elementary matrices are commutators. See \cite{Burger-Monod1,
Burger-Monod2} for a generalization to lattices in higher-rank
groups.
It is for this reason that we  work in dimension $2$.
It would be interesting to explore stable mixing
in higher dimensions.

We will prove Theorem~\ref{separation} with the use of elementary
hyperbolic
geometry. It seems to be likely that the separation property
stated in the theorem can be extended to general hyperbolic
groups -- the proof
could be probably extracted from a
recent work by Bestvina and Fujiwara \cite{BF}.

Dima Burago and David Kazhdan pointed out
to us that if we restrict the class of kicks
to sequences $\Phi$ taking a {\it finite}
number of values, then it is possible to show that the kicked system
\eqref{kicked system} is mixing for $t\gg 1$
under a more general condition then in
Theorem \ref{thm: mixing}: One need only assume
that the kicks do not exchange the stable and
unstable subspaces of the hyperbolic matrix $h$.
Their argument is of a dynamical nature, and  does not use
quasi-morphisms. No such argument is known
to us in the more general case of kicks
with bounded traces as in Theorem \ref{thm: mixing}.
In fact, stable mixing is sensitive
to the size of the kicks. For instance, stability disappears
when one allows certain unbounded sequence of kicks
(see Section 1.2 of [PR]).
It would be interesting
to determine the critical size of perturbations
for which mixing persists. We refer to \S 4 for further
discussion on this issue in terms of a special geometry
on $SL(2,\Z)$.

\noindent{\bf Acknowledgments:}
We thank Dima Burago, David Kazhdan
and Misha Sodin for several useful discussions.

\section{The separation theorem}

\subsection{Background on quasi-morphisms}
Let $G$ be a group.
Starting with any quasi-morphism $r'$ on $G$,
we get a homogeneous one $r$   by setting
$$
r(g) = \lim_{n\to \infty} \frac{r'(g^n)} n
$$
(the limit exists because $|r'(g^n)|$ is a subadditive sequence).
A homogeneous quasi-morphism is
automatically invariant under conjugation since for all $n$
$$
r(ghg^{-1}) = \frac{r((ghg^{-1})^n)}n  = \frac {r(gh^ng^{-1})}n =
\frac{r(g) + nr(h) - r(g) + O(1)}n
$$
and taking $n\to \infty$ we get $r(ghg^{-1}) =r(h)$.

\subsection{Quasi-morphisms for discrete subgroups of $PSL(2,\R)$}
Consider the action of
a discrete group $G\subset PSL(2,\R)$ on the
upper half-plane $\HH$.
To  construct  quasi-morphisms of $G$, start with
a smooth $G$-invariant one-form $\alpha$, 
such that $|d\alpha/\Omega|\leq C$ for some
$C>0$, where $\Omega = y^{-2}dx\wedge dy$ is the hyperbolic area form.
Given such $\alpha$ and a base-point $z\in \HH$,
set
$$
r_z(g) = \int_{\ell(z,gz)} \alpha
$$
where  $\ell(z,w)$ denotes the geodesic
segment between two points $z,w\in \HH$.
As is well known, $r_z$ is a quasi-morphism of $G$
(the reason that $dr_z$ is {\em bounded} has to do with the
fact that the area of a hyperbolic triangle is bounded - see
e.g. 
\cite[Lemma 3.2.B.]{kicks1}).
The quasi-morphism $r_z$ depends on $z$ only up to a bounded
quantity. Thus its homogenization
$$r(g):=\lim_n r_z(g^n)/n$$
is independent of the choice of base-point $z$.

\subsection{Proof of the separation theorem}
We will construct $\alpha$
so that for our fixed element $h$ and some point $z_0 \in \HH$,
$r_{z_0}(h^n) = n$. Hence for the homogenization $r$ one has
$r(h) = 1$.
Further, we make the support of $\alpha$ sufficiently
small so that  for each of the elements $g_j, \; j \geq 1$,
there is a base-point
$z_j \in \HH$ such that every geodesic segment $\ell(z_j,
g_j^n z_j)$ is disjoint from the support of $\alpha$. In that case
$r_{z_j}(g_j^n) = \int_{\ell(z_j,g_j^n z_j)}\alpha = 0$
and thus $r(g_j)=0$. This will conclude the proof of the Theorem.

We may assume that none of the elements $g_j$ are elliptic, since
these are annihilated by any homogeneous quasi-morphism. Thus we take
$g_1,\dots g_a$ to be hyperbolic, and $g_{a+1},\dots ,g_M$ parabolic.
 Let $L_1,\dots, L_a$ be the invariant geodesics of the hyperbolic
elements $g_1,\dots g_a$. Let $p_{a+1},\dots,p_M\in \R\cup\infty$ be
the cusps (fixed points) of the parabolic elements $g_{a+1},\dots
g_M$.

{\bf Case 1: $h$ primitive parabolic. } In this case $h$ is
not conjugate to its inverse already in $PSL(2,\R)$.
After conjugating the group $G$, we may assume
that $h(z)=z+1$. Since $G$ is discrete,
$\infty$ is not a fixed point of any hyperbolic
element of $G$. Therefore  if $R>1$ is sufficiently large,  the horoball
$B =\{z=x+iy:y>R\}$
is disjoint
from $\gamma B$ when $\gamma\in G$ is not a power of $h$.
Moreover, if $R>1$  is sufficiently large then $B$
contains none of the $G$-translates of
the geodesics $L_j$, $1\leq j\leq a$. Furthermore, since the
parabolic elements
$g_{a+1}, \dots g_M$ are not conjugate to a power of the primitive
element $h$, their fixed points $p_j$ are not $G$-translates of
$\infty$. Thus (further increasing $R$ if necessary)
there exist small horoballs
$B_j \subset \HH$ tangent to $\R$ at $p_j$ which are
disjoint from all $G$-translates of the horoball $B$.

We take a smooth cutoff function $u(y)$ which vanishes for $y<R$ and
is identically $1$ for $y>R+1$, and set $\alpha' =u(y)dx$.
Clearly $|d\alpha'/ \Omega|\leq C$.
Take $\alpha$ to be the $G$-periodization of $\alpha'$, that is
$$
\alpha = \sum_{\gamma\in \langle h \rangle \backslash
G}\gamma^*\alpha'.
$$
\noindent
This is a smooth one-form, supported in the union of $G$-translates
of the horoball $B$. Thus its support is disjoint from the
geodesics $L_j$ and the horoballs $B_j$.
For $j = 1,...,M$ pick a point $z_j$ which lies on $L_j$ for
$j \leq a$ and inside $B_j$ for $j > a$.
Thus $r_{z_j}(g_j)=0$
for all $1\leq j\leq M$.
To compute $r(h)$, take a point $z_0$
with ${\rm Im}(z_0)=R+1$. Then $h^n z_0=z_0+n$ and
$\int_{\ell(z_0,h^nz_0)}\alpha = n$, and so $r_{z_0}(h^n) = n$.
This proves the claim for parabolic $h$.

{\bf Case 2: $h$ primitive hyperbolic. }
Let $L$ be the invariant geodesic of $h$. Since $g_1,\dots, g_a $ are
not conjugate to a power of the primitive element $h$, their invariant
geodesics $L_j$ are not $G$-equivalent to $L$.
Let $D$ be a locally finite fundamental domain for $G$ (e.g. the
Dirichlet fundamental domain) whose intersection with $L$
contains a geodesic segment, say $I$.
There are only finitely many
$G$-translates of $L_j$ which meet $D$
\cite[Theorem 9.2.8(iii)]{Beardon}, and so the intersection of
$L\cap D$ with the $G$-translates of the $L_j$'s is a finite set of
points.
Shrinking $I$, we can assume that $I$
does not
meet any of these intersection points,
and moreover that a
small neighborhood $U$ of $I$ is disjoint
from the $G$-images of $L_j$.
Shrinking further $I$ and $U$ if necessary we can
achieve the following.
First, $\gamma U \cap U=\emptyset$ if $1\neq \gamma \in G$.
Further, since $h$ is not conjugate
to $h^{-1}$ in $G$, this implies that if
$\gamma\in G$, $\gamma\neq h^k$ for some $k\in \Z$ then $\gamma U$ is
bounded away from $L$ \cite[Lemma 3.2.D]{kicks1}.
Finally, we can choose small horoballs $B_j$ around the
cusps $p_j$ which are disjoint from $\cup_{\gamma\in G}\gamma U$.

Choose now a one-form $\alpha'$ in $\HH$, supported in $U$  so that
$\int_I\alpha'=1$. Then $|d\alpha'/\Omega|$ is bounded.
Let $\alpha =\sum_{\gamma\in G}\gamma^*\alpha'$ be
its $G$-periodization. Then the support of $\alpha$ lies in the union
of $G$-translates of $U$, and so is disjoint from the geodesics $L_j$
as well as from the horoballs $B_j$. For $j = 1,...,M$
pick a point $z_j$ which lies on $L_j$ when $j \leq a$ and
inside $B_j$ for $j > a$. Pick a point $z_0$ on $L$.
Clearly, $r_{z_j}(g_j) = 0$ for $j \geq 1$.
Finally note that
$\int_{\ell(z_0,h^nz_0)}\alpha = n\int_I \alpha = n$ and
so $r_{z_0}(h^n)=n$.
This concludes the claim for the hyperbolic case. \qed

\subsection{Proof of Proposition~\ref{prop: quasimorphism}}
Let $\pi:SL(2,\Z)\to PSL(2,\Z)$ be the projection.
Assume $h \in SL(2,\Z)$
is hyperbolic and not conjugate to its inverse in $SL(2,\Z)$.
Note that
$h$ is not conjugate to $-h^{-1}$ since $\tr (h)\neq 0$.
Let $\tilde r$ be a homogeneous quasi-morphism of $PSL(2,\Z)$ which is
$1$ on $\pi(h)$ and  zero on $\pi(\begin{pmatrix}1&1\\
&1\end{pmatrix})$ whose existence is guaranteed by
Theorem~\ref{separation}, and let $r=\tilde r\circ \pi$.
Since all primitive parabolic elements   are conjugate in $SL(2,\Z)$
to either $\begin{pmatrix}1&1\\&1\end{pmatrix}$ or its inverse, it
follows that $r$ vanishes on all parabolic elements of $SL(2,\Z)$. \qed

\subsection{When a matrix is conjugate to its inverse?}
Every symmetric matrix from $SL(2,\Z)$ is conjugate
to its inverse by an element
$\begin{pmatrix} 0 & 1\\ -1 & 0\end{pmatrix}.$
Further, one can prove the following statement, 
which is left as an exercise: For a hyperbolic matrix $h$, 
consider
the prime decomposition of $({\rm trace} (h))^2 - 4$.
If this decomposition contains a prime which equals
3 modulo 4 and enters in an odd power then $h$ is not conjugate
to $h^{-1}$. For instance,
$$ h =
\begin{pmatrix} 4 & 9\\ 7 & 16\end{pmatrix}$$
is not conjugate to its inverse (see also \cite{BR})
and thus is stably mixing.
There is an algorithm, based on the theory of Pell's equation, which
establishes when a given matrix from $SL(2,\Z)$ is conjugate to its
inverse, see \cite{Gottlieb, Long}.  

\section{Mixing and exponential decay of correlations}

\subsection{Decay of correlations}
In this section we prove Theorem~\ref{thm: mixing} as well as
a quantitative
version giving an exponential rate of decay of the
correlations functions
$C(F_1,F_2;f^{(n)})$
given by \eqref{correlation function}.
Mixing means that the correlation functions decay as $n\to \infty$ for
$L^2$ observables. To get a rate of decay, we have to restrict the
observables.
A function $F$ on the torus satisfies {\em H\"older's
condition} if for some $\gamma>0$ there is  a positive constant $c>0$
so that for all $x,y\in \TT$ we have
\begin{equation}\label{Holder condition}
|F(x)-F(y)|\leq c \norm{x-y}^\gamma \;.
\end{equation}
We will show  that our kicked systems
have exponential decay of correlations for H\"older
observables.
\begin{thm}\label{thm: decay of correlations}
If $h\in SL(2,\Z)$ is a hyperbolic and not conjugate to its inverse,
there is some $t_0$ so that for all $t>t_0$
and for all H\"older functions $F_1,F_2$ on $\TT$ the correlations
$C(F_1,F_2; f^{(n)}(t))$ decay exponentially:
$$
|C(F_1,F_2; f^{(n)}(t))|\leq c e^{-\gamma n}
$$
for some $c,\gamma>0$.
\end{thm}

In order to prove theorem~\ref{thm: decay of correlations}
we start with a general observation:
\begin{prop}\label{prop: expdecay}
Let $f^{(n)}\in SL(2,\Z)$ be a sequential system.
If there are $\alpha>0$, $\beta>1$ such that for every  integer
vector $0\neq v\in \Z^2$ and every $n\geq 1$ we have
$$
\norm{v  f^{(n)}}\geq \alpha \frac{\beta^n}{\norm{v}}
$$
then the system is mixing and moreover we have exponential decay of
correlations for all H\"older observables.
\end{prop}


\begin{proof}
It clearly suffices to deal with the case where $F_1=F_2$ and both
have mean zero. Let
\begin{equation}\label{choice of N}
N\leq \sqrt{\alpha \beta^n} \;.
\end{equation}
Take a function in $L^2(\TT)$, of mean zero,
and consider its Fourier expansion
$$
F(x) = \sum_{0\neq v\in \Z^2} a(v)
e^{2\pi i(v,x)}
 = F_N(x)+R_N(x)
$$
where $F_N(x) = \sum_{\norm{v}<N} a(v)
e^{2\pi i(v,x)}
$ is the partial sum over  all
frequencies of norm less than $N$.
Then the correlation function $C(F,F;f^{(n)})$ is a sum of terms
\begin{multline*}
C(F,F;f^{(n)}) =\int_{\TT} F_N (f^{(n)}x ) F_N(x)dx   \\
+\int_{\TT} F_N( f^{(n)}x) R_N(x)dx +
\int_{\TT} R_N( f^{(n)}x) F(x)dx.
\end{multline*}
By Cauchy-Schwartz, the second and third terms above are  bounded
respectively by
$$
\norm{F_N\circ f^{(n)}}_2 \norm{R_N}_2 = \norm{F_N}_2 \norm{R_N}_2
\leq \norm{F}_2 \norm{R_N}_2
$$
and by
$$
\norm{R_N\circ f^{(n)}}_2 \norm{F}_2 =\norm{R_N}_2\norm{F}_2 \;,
$$
where $||F||_2$ stands for the $L^2$-norm.
The first term equals
$$
\int_{\TT} F_N\circ f^{(n)}(x) F_N(x)dx  =
\sum_{\substack{\norm{vf^{(n)} }<N \\ \norm{v}<N}} a(-vf^{(n)})a(v).
$$
We claim that this is an empty sum, hence vanishes, for our choice of
$N$: Indeed, by our condition if $0<\norm{v}<N$ we have
$$
\norm{-vf^{(n)}}\geq \alpha \frac{\beta^n}{\norm{v}} >
\frac{\alpha\beta^n}N \geq N
$$
since $N\leq \sqrt{\alpha \beta^n}$. Thus we find that for
$N$ as above that
$$
|C(F,F;f^{(n)}| \leq 2 \norm{F}_2 \norm{R_N}_2 \;.
$$
Since we can choose $N\to\infty$ as $n\to \infty$ (subject to
\eqref{choice of N}) we get $\norm{R_N}_2\to 0$ and hence we have
mixing.

To prove {\em exponential} decay of correlations for H\"older
observables $F$,  recall
that if $F$ satisfies the H\"older condition (3)
of order $\gamma>0$ then for some constant
$c_F>0$
$$
\norm{R_N}_2 \leq c_F N^{-\gamma}
$$
(see e.g. \cite[Vol. I, Chapter II.4]{Zygmund} for the
proof in the case of functions of one variable )
and hence taking $N$ to be the smallest integer less than
$\sqrt{\alpha \beta^n}$ gives $\norm{R_N}_2 \leq
{\rm const}\; \beta^{-n\gamma/2}$ as
required.
\end{proof}

\subsection{Using quasi-morphisms}
For a quasi-morphism $r$ put
$$
\norm{dr}:= \sup_{g_1,g_2\in G} |r(g_1g_2)-r(g_1)-r(g_2)| \;.
$$
To show that our kicked systems satisfy the condition of
Proposition~\ref{prop: expdecay}, we use
\begin{prop}\label{prop: lower bound}
Let $r$ be a homogeneous quasi-morphism of $SL(2,\Z)$ which vanishes
on all parabolic elements.
Then there is $\alpha>0$ such that for every  integer
vector $0\neq v\in \Z^2$ and every $f \in SL(2,\Z)$, we have
$$
\norm{vA} \geq \alpha \frac {e^{|r(A)|/||dr|| }}{\norm{v}}
$$
\end{prop}

\noindent{\bf Proof:}
A vector $v\in \Z^2$ is {\em primitive} if
it is not a nontrivial multiple of
another integer vector.
It clearly suffices to prove the result for primitive vectors.
\begin{lem}\label{lem: gcd}
Every primitive vector $v\in \Z^2$ can be written as
$$
v= (0,1)\begin{pmatrix} 1 &k_N \\&1\end{pmatrix}
\begin{pmatrix} 1 & \\k_{N-1} &1\end{pmatrix} \dots
\begin{pmatrix} 1 &k_{2} \\&1\end{pmatrix}
\begin{pmatrix} 1 & \\k_{1} &1\end{pmatrix}
$$
with $N\leq \log_2\norm{v}+10$.
\end{lem}
\begin{proof}
This is the well known Euclidean algorithm:
Assume that $v=(p,q)$ with $|q|\geq |p|$.
Primitivity is equivalent to  $p,q$ being co-prime.
Then find an integer $k_1$ so
that $q=k_1 p+r_1$ with $|r_1|\leq |p|/2$. We then
arrive at a new vector
$v_2:=(p,r_1) = (p,q)\begin{pmatrix} 1 &-k_1 \\&1\end{pmatrix}$. Now
find an integer $k_2$ so that $p=k_2r_1 + r_2$ with $|r_2|\leq
|r_1|/2$ to get another vector $v_3:=(r_2,r_1) =
(p,r_1)\begin{pmatrix} 1 & \\ -k_2 &1\end{pmatrix} $ and so on. This
proceeds until we get to the point that we have computed the greatest
common divisor of $p$ and $q$ (which in our case equals $1$)
with at most $\log_2 |p|+1\leq \log_2 \norm{v}+1$  steps.
This gives us either the vector $(0,1)$, in which case we are done,
or the vector $(1,0)$. In the latter case just note that
$(0,1) =(1,0)\begin{pmatrix} 1 & 1\\ &1\end{pmatrix}
\begin{pmatrix} 1 & \\-1 &1\end{pmatrix}$.
\end{proof}


By lemma~\ref{lem: gcd} every primitive vector is of the form $w=(0,1)h$ with
$h$ a product of at most $\log_2\norm{w}+10$ elementary matrices.
The quasi-morphism $r$ vanishes on these matrices and so we have \\
$|r(h)|\leq (\log_2\norm{w}+10)\norm{dr}$.
Apply this reasoning to the vectors $v=(0,1)h_1$ and
$vf^{(n)}=(0,1)h_2$. Then
$|r(h_1)|\leq (\log_2\norm{v}+10)\norm{dr}$
and
$|r(h_2)|\leq (\log_2\norm{ vf^{(n)} }+10)\norm{dr}$.
The matrix $h_3:=h_1
f h_2^{-1}$ is parabolic, since  it fixes $(0,1)$ and so
\begin{equation*}
\begin{split}
|r(f)|=|r(h_1^{-1}h_2h_3)| &\leq |r(h_1)|+|r(h_2)|+ 2\norm{dr}\\
&\leq (\log_2\norm{v} + \log_2\norm{vf} + 22) \norm{dr}
\end{split}
\end{equation*}
Therefore
$$
\norm{vf} \geq 2^{-22} \frac{ 2^{ |r(f) |/ ||dr||}}
{ \norm{v}}
$$ as required.
\qed

\subsection{Proof of Theorem~\ref{thm: decay of correlations}}
By Proposition~\ref{prop: quasimorphism} if $h$ is not
conjugate to its inverse then there is a homogeneous quasi-morphism
$r$ of $SL(2,\Z)$ which vanishes on parabolic elements and for which
$r(h)=1$. We claim that for any sequence
of kicks $\phi_i \in G$
with bounded traces one has $|r(\phi_1)| \leq c$ for some $c > 0$.
Indeed, $r$ vanishes on elliptic and parabolic elements, while
hyperbolic elements whose traces are bounded represent a {\it finite}
number of conjugacy classes in $G$. The claim follows.

Further,
$$
r(f^{(n)}(t)) = r(\prod_{i=1}^n h^t\phi_i)  =
\sum_{i=1}^n r(h^t)+r(\phi_i) + O_r(1) =
nt +  O_{r,\Phi}(n)
$$
and so if $t$ is big enough we have  $|r(f^{(n)}(t))|\geq n$.
By proposition~\ref{prop: lower bound} it follows that for all nonzero
integer vectors $v$ we have
$$
\norm{vf^{(n)}(t)}\geq \alpha \frac{\beta^n}{\norm{v}}
$$
for some $\beta>1$, $\alpha>0$
and thus Proposition~\ref{prop: expdecay} concludes the proof.
\qed

\section{Stable mixing from the geometric viewpoint}

In this section we present a generalization of
the stable mixing phenomenon
described above. Define a biinvariant
metric $\rho$ on the group
$G = SL(2,\Z)$ as follows.
Every element $g \in G$ can be decomposed
as $g = a_1...a_d$, where every
$a_i$ is either elliptic, or parabolic, or a
simple commutator of the form $aba^{-1}b^{-1}$, $a,b \in G$.
Put $\rho (1,g) = \inf d$,
where the infimum is taken over all such presentations, and set
$\rho(f,g) = \rho (1, fg^{-1})$ for all $f,g \in G$.

Denote by $G^{\infty}$ the set of all infinite
sequences $\{g^{(k)}\},\; k \in \N$
of "moderate growth", that is of those
which satisfy $\sup_{k \in \N} \rho(1,g^{(k)})/k < \infty.$
For instance, for every $h \in G$ the cyclic semigroup
$\{h^k\},\; k \in \N$
lies in $G^{\infty}$. Define a metric on $G^{\infty}$
by
$$\bar \rho
(\{f^{(k)}\}, \{g^{(k)}\}) = \sup_{k \in \N} \rho (f^{(k)},g^{(k)}).$$

\begin{thm}\label{thm: geom} Let $h \in G$ be a hyperbolic element
which is not conjugate to its inverse in $G$. Then
there exists $\eps > 0$
such that the ball of radius $\eps$ centered at
$\{h^k\} \in G^{\infty}$
consists of mixing sequences with exponential decay of correlations
on H\"older observables.
\end{thm}

\begin{proof} One can easily check that
homogeneous quasi-morphisms on $G$
which vanish on parabolics
are Lipschitz functions with respect to metric $\rho$. Let $r$ be
such a quasi-morphism  with $r(h) = 1$. Take $\eps > 0$
sufficiently small  and consider any
sequence $ f = \{f^{(k)}\} \in G^{\infty}$
with $\bar \rho(\{f^{(k)}\},\{h^k\}) < \eps$. It follows
that $|r(f^{(k)})| \geq ck$ for some $c > 0$ and all $k \in \N$.
Applying
Propositions 7 and 6, we get that $f$ is mixing
with exponential decay of
correlations on H\"older observables.
\end{proof}

\medskip
\noindent
{\bf Remark.}
Throughout the paper we worked with hyperbolic elements $h$ which
are not conjugate to $h^{-1}$ in $G$. This assumption can be translated
in the geometric language as follows.
We say that a cyclic subgroup $\{h^k\}, \; k \in \Z$
has {\it linear growth} if $\rho(1,h^k) \geq c|k|$ for some $c > 0$.
It follows from Proposition 3 that the linear growth is equivalent to
the fact that $h$ is a hyperbolic element which
is not conjugate to its inverse in $G$. Otherwise cyclic
subgroup generated by $h$ remains a bounded distance from the  identity.

\medskip

To include kicked systems into the geometric framework introduced
above, we start with the following observation. Denote by $\eps(h)$
the maximal value of $\eps$ supplied by Theorem \ref{thm: geom}.
Then for $t \in \N$ one has $\eps(h^t) \geq t\eps(h)$.
Assume now that $h$ is not conjugate to its inverse in $G$,
take a
sequence of kicks $\Phi = \{\phi_i\}$
with bounded traces, and consider the kicked system
$f^{(n)}(t) = \phi_nh^t...\phi_1h^t$.
Note that
$$h^{nt} = f^{(n)}(t) \circ \prod_{j=1}^n h^{-tj}\phi_j^{-1}h^{tj}.$$
Therefore there exists $c > 0$ which depends only on $\Phi$ such that
for every $t \in \N$
$$\bar \rho (\{h^{nt}\},\{f^{(n)}(t)\}) \leq c.$$
Hence the kicked system $f^{(n)}(t)$ is mixing provided
$c < t\eps(h) \leq \eps(h^t)$.

\section{Positive Lyapunov exponent}

We say that a sequential system $\{f^{(n)}\}$ has {\it positive Lyapunov
exponent} if the sequence ${\rm trace}(f^{(n)})$ grows exponentially
with $n$. It turns out that
kicked system (1) has positive Lyapunov exponent
provided $h$ is not conjugate to its inverse
and $t$ is large enough. Indeed, choose a homogeneous quasi-morphism $r$
which satisfies $r(h)=1$ and vanishes
on parabolics, and apply the following
estimate.

\begin{thm}
$$|{\rm trace} (f)| \geq (2{\sqrt 5})^{|r(f)|/||dr||},$$
for every hyperbolic $f \in SL(2,\Z)$ and any homogeneous quasi-morphism
$r$ which vanishes on parabolics.
\end{thm}

\begin{proof} The proof is divided into 3 steps.

1) We claim that every hyperbolic matrix
$f = \begin{pmatrix} a & b\\c & d\end{pmatrix}$
is conjugate in $SL(2,\Z)$
to a matrix
$g = \begin{pmatrix} a_1 & b_1\\c_1 &d_1\end{pmatrix}$
with $|c_1| \leq |{\rm trace} (f)|/\sqrt 5$.
\noindent
Indeed, assume that $g = hfh^{-1}$, where
$h = \begin{pmatrix} z & t\\x & y\end{pmatrix} \in
SL(2,\Z).$ One calculates that
$c_1 = Q(x,y)$, where $Q(x,y) = cx^2 + (a-d)xy - by^2$.
The discriminant $\Delta$ of the quadratic form $Q$ equals
${\rm trace}^2 (f) - 4$. Since $f$ is hyperbolic, $\sqrt \Delta$ is
irrational. In this case (see \cite{Cassels}) there exists a primitive
vector $(x,y) \in \Z^2$ with
$|Q(x,y)| \leq {\sqrt \Delta}/{\sqrt 5} \leq |{\rm trace} (f)| /\sqrt 5.$
Hence we obtain $h \in SL(2,\Z)$ as required.

2) As a consequence we get that every hyperbolic matrix $f \in SL(2,\Z)$
decomposes as $f = f'h$, where
$|{\rm trace} (f')| \leq |{\rm trace} (f)|/(2{\sqrt 5})$
and $h$ is parabolic.
 Indeed, in view of step
1 we can assume without loss of generality that
$f = \begin{pmatrix} a & b\\c & d\end{pmatrix}$
with $|c| \leq |{\rm trace} (f)|/\sqrt 5$.
Write
$f' =f \begin{pmatrix} 1 & k\\ & 1\end{pmatrix},$
so ${\rm trace} (f') = {\rm trace} (f) + ck$. One can always choose
$k \in \Z$ so that
$|{\rm trace} (f')| \leq |c|/2 \leq |{\rm trace} (f)|/(2{\sqrt 5})$.

3)
For an integer $s \geq 3$ put $u(s)=
\max |r(g)|$, where the maximum is taken over all hyperbolic matrices $g$
whose trace lies in $[-s;s]$. Since a quasi-morphism is a class function on the
group, and the number of hyperbolic
conjugacy classes with given trace is finite, the
function $u$ is well defined. Step 2 yields inequality
$u(s) \leq u([s/(2\sqrt 5)]) + ||dr||$ for all $s \geq 3$,
where we set $u(s) = 0$ for $s \leq 2$, and brackets stand for the
integral part of a  real number.
Arguing by induction we get that
$u(s) \leq ||dr|| \log_{2\sqrt 5} s$ for all $s \geq 3$.
This completes the proof.
\end{proof}


\end{document}